\documentclass[11pt]{amsart}

\usepackage[a4paper,margin=3.5cm, centering]{geometry}
\usepackage[T1]{fontenc}
\usepackage[euler-digits,euler-hat-accent]{eulervm}
\usepackage[]{ccicons}

\title[A bigroupoid's topology]{A bigroupoid's topology\\
{\tiny (or, Topologising the homotopy bigroupoid of a space)}}

\author{David Michael Roberts}
\thanks{This work was undertaken while the author was supported by an Australian Postgraduate Award and by Australian Research Council grant number DP120100106. Further financial support was provided by Mrs R. \\
\cczero\ This article is released under a CC0 license, \url{https://creativecommons.org/publicdomain/zero/1.0/}}
\email{david.roberts@adelaide.edu.au}
\address{School of Mathematical Sciences, University of Adelaide, 5005 Australia}

\usepackage{color}
\usepackage{subfig}
\usepackage{graphicx}

\usepackage{amsmath}
\usepackage{amssymb}
\usepackage{amsthm} 
\usepackage{enumerate}

\usepackage[parfill]{parskip}

\newtheorem{theorem}{Theorem}[section]
\newtheorem{lemma}[theorem]{Lemma}
\newtheorem{proposition}[theorem]{Proposition}
\newtheorem{corollary}[theorem]{Corollary}
\theoremstyle{definition}
\newtheorem{definition}[theorem]{Definition}
\newtheorem{example}[theorem]{Example}
\newtheorem{remark}[theorem]{Remark}

\usepackage{hyperref}

\usepackage{ifpdf}
\ifpdf
\usepackage[all,pdf]{xy} 
\else
\input xy
\xyoption{all}
\xyoption{2cell}
\xyoption{v2}
\fi
\UseComputerModernTips

\newcommand{\st}{\rightrightarrows}
\newcommand{\into}{\hookrightarrow}
\newcommand{\cHom}[1]{\mathcal{H}om_{#1}}
\newcommand{\W}{\mathcal{W}}
\newcommand{\Y}{\Upsilon}

\DeclareMathOperator{\disc}{disc}
\DeclareMathOperator{\id}{id}
\DeclareMathOperator{\pr}{pr}
\DeclareMathOperator{\ev}{ev}
\DeclareMathOperator{\im}{im}
\DeclareMathOperator{\Aut}{Aut}
\DeclareMathOperator{\Obj}{Obj}
\DeclareMathOperator{\Mor}{Mor}
\DeclareMathOperator{\Top}{\textbf{Top}}
\DeclareMathOperator{\Gpd}{\textbf{Gpd}}
\DeclareMathOperator{\Bigpd}{\textbf{Bigpd}}

\begin{document}

\begin{abstract}

  The fundamental bigroupoid of a topological space is one way of capturing its homotopy 2-type. 
  When the space is semilocally 2-connected, one can lift the construction to a bigroupoid internal to the category of topological spaces, as Brown and Danesh-Naruie lifted the fundamental groupoid to a topological groupoid.
  For locally relatively contractible spaces the resulting topological bigroupoid is \emph{locally trivial} in a way analogous to the case of the topologised fundamental groupoid.


\end{abstract}
\maketitle

\section{Introduction}
\label{sec:intro}

One of the standard examples of a groupoid is the fundamental groupoid $\Pi_1(X)$ of a topological space $X$, generalising the fundamental group $\pi_1(X,x)$ at a basepoint $x$ to consider `all basepoints at once'. 
In \cite{Brown-Danesh-Naruie_75} Brown and Danesh-Naruie showed that, under a mild assumption, the fundamental groupoid can be given the structure of a \emph{topological} groupoid. 
That is, the sets of objects and arrows---points in the space and homotopy classes of paths, respectively---can be given topologies such that all the maps that make up the groupoid (source, target, composition etc) are continuous. 

The mild assumption mentioned in the previous paragraph is exactly that which guarantees the existence of a universal covering space; a seemingly little-known fact is that said covering space can be constructed directly from the topologised fundamental groupoid as given in \cite{Brown-Danesh-Naruie_75}.
Moreover, this construction is formally analogous to the construction of the first stage of the Whitehead tower of a topological space.

Drawing inspiration from the celebrated Homotopy Hypothesis linking higher groupoids and homotopy types, we see that to extend these constructions to dimension 2 we need to consider some form of 2-dimensional groupoid. 
While there are several different algebraic models that completely capture the homotopy 2-type of a space, such as crossed modules (Whitehead, 1940s) and double groupoids (Brown--Higgins 1970s), here we choose to consider \emph{bigroupoids}; Stevenson \cite{Danny_phd} and Hardie--Kamps--Kieboom \cite{HKK_01} constructed a fundamental bigroupoid $\Pi_2(X)$ of a space $X$. 
The idea of such an object, albeit in the fully general case of weak $n$-groupoids representing arbitrary homotopy $n$-types, seems to go back to Grothendieck's 1975 letters to Breen \cite{Grothendieck_75}.

The idea of a bigroupoid is illustrated nicely by considering this special case.
Firstly, bigroupoids have object and arrows, as groupoids do, but also \emph{2-arrows}, which are arrows between arrows.
Objects of $\Pi_2(X)$ are points in $X$ and arrows are paths $I = [0,1] \to X$.
Paths can be composed, but since at this point there is no quotient by the relation of homotopy, such composition is not associative.
Similar issues arise when composing by constant paths, or reverse paths, representing identity arrows and inverses repectively.
This is where the 2-arrows come in: 2-arrows in $\Pi_2(X)$ are homotopy classes of homotopies of paths. 
Or, equivalently, homotopy classes of \emph{bigons}, which are certain maps $I^2 \to X$.
As maps $I^2 \to X$ support pasting in two directions, we get the \emph{horizontal} composition of bigons end-to-end (inducing composition on their boundary paths) and the \emph{vertical} composition of bigons pasted along one of their boundary paths.

This article will give, under a mild local condition, topologies for all the sets involved in $\Pi_2(X)$---points in $X$, paths in $X$ and homotopy classes of bigons---such that every operation in the bigroupoid structure is continuous. 
One of the reasons that a strict model is not chosen is that they do not seem well-adapted to the application that motivated the present author, namely constructing geometrically and in a smooth fashion the \emph{second} stage in the Whitehead tower of a manifold.
Double groupoids and crossed modules over groupoids, both championed by Ronnie Brown, seem to work best in the context of computing with topologically discrete algebraic structures (see however the concluding remarks in section \ref{sec:local_triv}). 
Likewise the homotopy \emph{2-groupoid} of a Hausdorff space given in \cite{HKK_00} uses thin homotopy classes which does not lead to a well-behaved space of arrows\footnote{Even worse, in the smooth setting, one does not even have a half-decent manifold structure on the set of thin homotopy classes of loops, see \cite{Stacey_10}.}.

The approach of the paper is that one can in fact take the given topologies on the sets of objects and arrows of $\Pi_2(X)$, namely the topology on $X$ and the compact-open topology on $X^I$;
the main novelty is to define a very particular basis for the topology on the set of homotopy classes of bigons so that one can prove the required continuity of structure maps involving 2-arrows.
This uses in an essential way the local assumptions on $X$.
The paper finishes by showing that $\Pi_2(X)$, with the topology we define, satisfies analogues of the local triviality\footnote{Local triviality of topological groupoids is a condition, introduced by Ehresmann \cite{Ehresmann_59}, that relates them with locally trivial principal bundles.} and discreteness properties that the topological groupoid $\Pi_1(X)$ has.

Extending these results further up the ladder of higher groupoids needs to take a different approach, because even weak 3-groupoids---\emph{strict} 3-groupoids are known to be insufficient---are quite complicated. After that, the explicit algebraic definitions are no longer practical if one wants to capture the full homotopy type.
One could consider however other models for higher groupoids, such as operadic definitions of weak $n$-groupoids; the approach of Trimble \cite{Trimble_99} seems like it may be appropriate, given the approach of the sequel \cite{Roberts_15b} to this paper.
The analogue of the results in the current paper would be that, under suitable local connectivity assumptions, the algebras for the operads involved in the definitions would be \emph{topological}, i.e.~algebras in the category of spaces rather than in the category of sets.
Alternatively one might use Kan complexes with certain unique filler conditions and then consider internal Kan complexes in $\Top$, or even simplicial sheaves on $\Top$, as models for higher topological groupoids.

In \cite{Bakovic_phd} Bakovic gives a recipe, partly building on \cite{Roberts-Schreiber_08}, for taking an internal bigroupoid (for instance in topological spaces) and giving a principal 2-bundle. 
Topological bigroupoids with non-discrete object space do not seem to be very common, so this paper gives at least one family of examples for Bakovic's general machinery.
In fact the resulting principal 2-bundle is the desired second stage in the Whitehead tower, as was constructed in the author's thesis \cite{Roberts_phd}.

Thanks are due to several anonymous referees who helped beat this article into shape over several iterations, and to Tim Porter for both inviting its submission to this volume and his subsequent patience.
Thanks also to Ronnie Brown, whose lovely book \cite{RBrown_groupoids} on groupoids and topology was influential in my thesis work (of which this paper formed a small part) in ways that are not apparent to the casual observer: Happy Birthday Ronnie!

\section{Topological groupoids and bigroupoids}
\label{sec:top_gpds_bigpds}

Recall that a topological groupoid is a groupoid with a space of objects and a space of arrows such that all the structure maps are continuous. 
Functors $\Gamma\to \Delta$ between topological groupoids $\Gamma$ and $\Delta$ consist of \emph{continuous} maps $\Gamma_0 \to \Delta_0$ and $\Gamma_1 \to \Delta_1$ commuting with all the groupoid structure.
The reason that we do not use the term `continuous functor' here is that this has a separate meaning for functors unrelated to topology.
The category of topological groupoids will be denoted by $\Gpd(\Top)$.

Recall that there is a full inclusion $\disc\colon \Top \to \Gpd(\Top)$, sending a topological space $X$ to the topological groupoid $\disc(X)$, with arrows and objects both given by $X$ with all structure maps the identity. 
All `spaces' will be topological spaces in what follows, unless otherwise specified.

To describe the topological fundamental bigroupoid $\Pi_2(X)$ of a space $X$, we first need to define topological bigroupoids. 
Such a thing may be defined using the full diagrammatic definition of an \emph{internal} bicategory in $\Top$ as in B\'enabou's \cite{Benabou}, which gives all the structure maps and spaces explicitly together with many commuting diagrams.
We will adopt instead a more compact approach.
For those familiar with such things, the definition below is the internal analogue of weak enrichment in groupoids.
For the uninitiated, one can think of the definition of a (topological) bigroupoid as being a generalisation of the following reworking of the definition of groupoid.

An ordinary groupoid $\Gamma$ is given by a set $\Gamma_0$ together with a family of sets $\Gamma(a,b)$, the set of arrows from the object $a$ to the object $b$. 
One should think of this as a set $\Gamma_1$ parameterised by $\Gamma_0 \times \Gamma_0$, or in other words, a set \emph{over} $\Gamma_0\times \Gamma_0$, written $(s,t)\colon \Gamma_1 \to \Gamma_0 \times \Gamma_0$. Composition is given by a function
\[
    \Gamma_1 \times_{\Gamma_0} \Gamma_1 \to \Gamma_1
\]
respecting the maps down to $\Gamma_0\times \Gamma_0$. Here the pullback is $\{(f,g)\in \Gamma_1\times \Gamma_1 \mid t(f) = s(g)\}$, considered as a set over $\Gamma_0\times \Gamma_0$ via $(f,g) \mapsto (s(f),t(g))$. 
Associativity can be enforced by asking that a certain diagram in sets over $\Gamma_0\times \Gamma_0$ commutes.
Likewise, inversion in the groupoid is an endomorphism of $\Gamma_0$ covering the swap map on $\Gamma_0\times \Gamma_0$, and if $\Gamma_0$ is considered as a set over $\Gamma_0\times \Gamma_0$ by the diagonal map, then the function assinging identity arrows is the map $\Gamma_0 \to \Gamma_1$.

Moving to bigroupoids, the hom-sets are replaced by hom-\emph{groupoids}, as can be seen by considering the case of $\Pi_2(X)$. 
For two fixed objects $x$ and $y$---points in $X$---we have a set of paths from $x$ to $y$, and a set of (homotopy classes of) bigons with vertices $x$ and $y$, and such bigons can be pasted vertically along a common edge.
This, together with degenerate bigons and reversal of orientation gives a groupoid.
Now, allowing $x$ and $y$ to vary we see that what we have is a family of groupoids parameterised by $X\times X$, or, in other words, a groupoid equipped with a functor to $\disc(X\times X)$.
Horizontal composition can then be encoded by a functor, and this composition is now \emph{not} associative. 
The commuting diagram of functions between sets that encodes associativity is now a diagram of functors between groupoids and only commutes up to a natural isomorphism, which of course needs to satisfy coherence conditions.
A generalisation of this approach was used by Trimble \cite{Trimble_99}, for instance, to define a general notion of weak higher groupoid.

The definition of \emph{topological} bigroupoid $B$ takes this idea of a family of hom-groupoids and replaces it by a continuous family of topological groupoids over the space $B_0\times B_0$, or in other words, a functor $(S,T)\colon \cHom{B} \to \disc(B_0\times B_0)$ between topological groupoids.
This definition can be unpacked to recover the standard definition of a bigroupoid, but would take up a fair amount of space.\footnote{For the sake of consiseness, any pullbacks or iterated pullbacks  over the space $B_0$ will follow the following convention: letting $H = \cHom{B}$ or $\Obj(\cHom{B})$, pullbacks of the form $(-) \times_{B_0} H$ use the functor $S\colon \cHom{B}\to \disc(B_0)$ or its object component, and pullbacks of the form $H\times_{B_0} (-) $ use the map $T\colon \cHom{B} \to \disc(B_0)$ or its object component.}
In the following definition $B_0$ is a stand-in for $\disc(B_0)$ when necessary to save space.

\begin{definition}\label{defn:top_bigpd}

  A \emph{topological bigroupoid} $B$ is a topological space $B_0$ (the \emph{space of objects}) and a topological groupoid $\cHom{B}$ (the \emph{hom-groupoid}, with source and target maps denoted $s_1$, $t_1$ respectively) equipped with a functor $(S,T)\colon \cHom{B} \to \disc(B_0\times B_0)$, together with:

  \begin{enumerate}[--]
    \item
      functors
      \begin{align*}
        \bullet\colon & \cHom{B} \times_{B_0} \cHom{B} \to \cHom{B}\\
        I\colon &\disc(B_0) \to \cHom{B}
      \end{align*}
      (\emph{composition} and \emph{identity}, respectively) over $B_0\times B_0$ and a functor
      \[
        \overline{(\cdot)}\colon \cHom{B} \to \cHom{B}
      \]
      (\emph{inverse}) covering the swap map for $B_0\times B_0$;
    \item
      letting $B_1 :=\Obj(\cHom{B})$ and $B_2 := \Mor(\cHom{B})$, continuous maps
      \begin{equation} \label{structure_maps}
        \begin{array}{l}
          a\colon B_1 \times_{B_0}  B_1\times_{B_0}  B_1 \to B_2 \\
          r\colon B_1 \to B_2\\
          l\colon B_1 \to B_2\\
          e\colon B_1 \to B_2\\
          i\colon B_1 \to B_2
        \end{array}
      \end{equation}
      that are the component maps of natural isomorphisms
      \[
        \xymatrix{
          \cHom{B} \times_{B_0} \cHom{B}\times_{B_0} \cHom{B}
          \ar[rr]^-{\id\times \bullet} \ar[dd]_{\bullet\times \id}
          &&
          \cHom{B} \times_{B_0} \cHom{B}
          \ar[dd]_{\ }="t"^{\bullet}\\
          \\
          \cHom{B} \times_{B_0} \cHom{B}
          \ar[rr]^(.8){\ }="s"_-{\bullet}
          &&
          \cHom{B}
          \ar@{=>}"s";"t"^a
        }
      \]\\
      \[
        \xymatrix@C-1pc{
          \cHom{B} \times_{B_0} \disc(B_0)
          \ar[rr]^-{\id\times I} \ar[rrdd]^{\ }="t1"_\simeq &&
          \cHom{B} \times_{B_0} \cHom{B}
          \ar[dd]^{\bullet}="s2"_{\ }="s1" && \disc(B_0) \times_{B_0} \cHom{B}
          \ar[ll]_-{I\times\id} \ar[ddll]_{\ }="t2"^\simeq\\
          \\
          && \cHom{B}
          \ar@{=>}"s1";"t1"_r
          \ar@{=>}"s2";"t2"^l
        }
      \]\\
      \[
        \xymatrix{
          \cHom{B} \ar[rr]^-{(\overline{(\cdot)},\id)} \ar[dd]_S &&
          \cHom{B} \times_{B_0} \cHom{B}
          \ar[dd]_{\ }="s1"^{\bullet}="t2" &&
          \cHom{B} \ar[ll]_-{(\id,\overline{(\cdot)})} \ar[dd]^T\\
          \\
          \disc(B_0)\ar[rr]^{\ }="t1"_I && \cHom{B} && \disc(B_0)\ar[ll]_{\ }="s2"^I
          \ar@{=>}"s1";"t1"_e
          \ar@{=>}"s2";"t2"_i
        }
      \]
  \end{enumerate}

  These are required to satisfy the usual coherence diagrams, for which the reader can refer to \cite[definitions 8.1, 8.2]{Danny_phd} (for instance).

\end{definition}

We can also define \emph{strict} 2-functors between bigroupoids. 
There is of course a notion of weak 2-functor between bigroupoids (called in \cite{HKK_01} a `pseudo functor'), but our functoriality results give strict 2-functors, so this is all that is needed here.

\begin{definition}\label{def:strict_2-functor}

  Let $B$ and $D$ be a pair of bigroupoids. A \emph{strict 2-functor} $F\colon B\to D$ is given by a continuous map $F_0\colon B_0 \to D_0$ and a functor $\mathcal{F}\colon \cHom{B}\to \cHom{D}$ of topological groupoids covering the induced map $B_0\times B_0\to D_0\times D_0$.
  This map is required to commute with the functors $\bullet$, $I$ and $\overline{(\cdot)}$ on each side, as well as respect the natural transformations $a,r,l,e$ and $i$.

\end{definition}

If we ignore the topology, Definition~\ref{defn:top_bigpd} is equivalent to the usual definition of a bigroupoid (for instance \cite[Definition 1.3]{HKK_01}), by considering individual hom-groupoids (that is, the fibres of $(S,T)$) and the induced functors thereon.
Compare the treatment in \cite[\S 1.5]{Leinster_03}, which defines bicategories as weakly enriched categories.

We define the (1-)category of topological bigroupoids and continuous strict 2-functors and denote it by $\Bigpd(\Top)$.

\section{The topological fundamental bigroupoid of a space}
\label{sec:top_Pi2}

A full definition of the fundamental bigroupoid $\Pi_2(X)$ can be found in \cite[Example~8.1]{Danny_phd} or \cite[\S2]{HKK_01}. 
We shall define it along the lines of Definition~\ref{defn:top_bigpd} as follows, noting that once the definitions are matched up, one gets an identical bigroupoid.
To distinguish the bare bigroupoid with no topology from the topologised version given below, we shall denote the former by $\Pi_2^\delta(X)$.

We need to define the `mild local condition' mentioned in the introduction:

\begin{definition}
\label{def:sl2c}

	A topological space is \emph{semilocally 2-connected} if it has a neighbourhood basis consisting of simply-connected sets $U$ with the inclusion $U \into X$ inducing the zero map $\pi_2(U) \to \pi_2(X)$, for any choice of basepoint.

\end{definition}

For instance, any locally contractible space like a manifold or CW-complex is semilocally 2-connected. 
Conversely, any semilocally 2-connected space is semilocally simply-connected.

Let $\Top_{sl2c}$ denote the full subcategory of $\Top$ on the semilocally 2-connected spaces.
We now make for the rest of the paper the assumption that $X$ is semilocally 2-connected, and only define the topological fundamental bigroupoid for such spaces. 
To start with, the space of objects $\Pi_2(X)_0$ is just the space $X$.
We need to then define the hom-groupoid $\cHom{\Pi_2(X)}$, as a groupoid over $\disc(X\times X)$.
It is built as follows:

\begin{enumerate}[--]

\item The space of objects of $\cHom{\Pi_2(X)}$ is $X^I$, the path space of $X$ with the compact-open topology.

\item
We define a \emph{bigon} to be a map $f\colon I^2 \to X$ that is constant on $\{\epsilon\}\times [0,1] \into I^2$ for $\epsilon=0,1$. 
Homotopy of bigons will always be relative to the boundary, so that homotopic bigons have equal boundaries.

\item The (underlying set of the) space $\Pi_2(X)_2$ of arrows of the hom-groupoid is the set of homotopy classes of bigons. 
Such homotopy classes will be referred to as \emph{2-tracks}, and written as $[f]$, for $f$ a representing bigon.
The source path is the restriction of the bigon to $[0,1]\times \{0\}$, and the target is the restriction to $[0,1]\times\{1\}$.
The topology on this set will be defined below in Subsection~\ref{subsec:topology_on_2-arrows}.

\item Composition in the hom-groupoid is by pasting 2-tracks in the direction of the second co\"ordinate; the identity 2-arrow is represented by the constant bigon on a path; inverses are given by precomposing with $(s,t) \mapsto (s,1-t)$, reversing the direction of a representing bigon (see Subsection~\ref{subsec:top_bigpd}). Denote this composition operation by $+$ and inversion with respect to it by $-(\cdot)$.

\item On objects the functor $(S,T)\colon \cHom{\Pi_2(X)} \to \disc(X\times X)$ is evaluation at the endpoints, which is continuous, and on arrows it is the composite $\Pi_2(X)_2 \xrightarrow{s_1} X^I \xrightarrow{(\ev_0,\ev_1)} X\times X$, sending $[f]\mapsto \left(f(0,0),f(0,1)\right)$. 
Hence to prove continuity we only need to show $\Pi_2(X)_2 \to X^I$ is continuous (see Subsection~\ref{subsec:hom-gpd_source_target_cts}).

\end{enumerate}

The next part of the definition is the composition, identity and horizontal inverse functors. 
The second of these is easy: it is simply the constant-path map $X\to X^I$, which is continuous.
The horizontal inverse functor $\cHom{\Pi_2(X)} \to \cHom{\Pi_2(X)}$ is, on objects, the reverse path map $X^I \to X^I$, which is continuous and manifestly covers the swap map on $X\times X$.
On morphisms this sends a 2-track $[f]$ to the homotopy class of the bigon $(s,t) \mapsto f(1-s,t)$.
Given the definition of $(S,T)$ this also covers the swap map. 

The horizontal composition functor on objects is concatenation of paths---again continuous---and on morphisms it is given by concatenating representative bigons in the direction of the first co\"ordinate.
Horizontal composition will be denoted by $\bullet$, and will be shown to be continuous on 2-tracks in Theorem~\ref{thm:we_have_a_top_bigpd} below.
The component maps (\ref{structure_maps}) of the natural isomorphisms in Definition~\ref{defn:top_bigpd} are given in detail in \cite[Example~8.1]{Danny_phd}, but can be reconstructed from any book that gives a definition of the fundamental group; for instance the associator $a$ is the 2-track with representative bigon the usual homotopy that encodes associativity of $\pi_1$.

The rest of this section is thus devoted to defining the topology on $\Pi_2(X)_2$ (Subsection~\ref{subsec:topology_on_2-arrows}), that the hom-groupoid $\cHom{\Pi_2(X)}$ is a \emph{topological} groupoid (Subsections~\ref{subsec:hom-gpd_source_target_cts} and \ref{subsec:hom_gpd_topological}) with a continuous functor to $\disc(X\times X)$ and that $\Pi_2(X)$ is a topological bigroupoid (Subsection~\ref{subsec:top_bigpd}).

\subsection{Topology on the set of 2-tracks}
\label{subsec:topology_on_2-arrows}

To describe the topology on the set of 2-tracks we will use a particular class of basic open sets of $X^I$ (and $X^{S^1}$) as follows.

\begin{definition}
\label{def:basis_for_comp_open_topology}
  Let $\gamma \in X^I$, and let $\mathfrak{p} = \{0<a_1<\ldots<a_{2n}<1\}$ be a partition of the unit interval.
  Also, let $\W=\{W_i\}_{i=1}^{2n+1}$ be a collection of basic open neighbourhoods in $X$ such that
  \begin{enumerate}[--]
    \item
      $W_{2i} \subset W_{2i-1}\cap W_{2i+1}$ for $i=1,\ldots,n$, and

    \item
      $\gamma([a_{i-1},a_i]) \subset W_i$ (by convention, let us take $a_0 = 0,\ a_{2n+1} = 1$).

  \end{enumerate}
  Define the set $N(\mathfrak{p},\W)\subset X^I$ to consist of those paths $\eta$ such that $\eta([a_{i-1},a_i]) \subset W_i$ for $i=1,\ldots,2n+1$.

  Similarly, if $\gamma \in LX = X^{S^1}$, one can further ask, given a collection $\W = \{W_i\}_{i=0}^{2n+1}$ of basic open neighbourhoods of $X$, that $W_0 \subset W_{2n+1} \cap W_1$. Define the set $N^o(\mathfrak{p},\W) \subset LX$ to consist of those loops $\omega$ such that $\omega([a_{i-1},a_i])\subset W_i$, where now $i$ is considered modulo $2n+2$.

\end{definition}

These two families of sets are shown to be a system of basic open neighbourhoods for the compact-open topology on $X^I$ and $LX$ in \cite{Roberts_phd}, Propositions 5.10 and 5.15 respectively.
When $X$ is semilocally 2-connected the sets $N(\mathfrak{p},\W)$ and $N^o(\mathfrak{p},\W)$ are relatively 1-connected: they are path connected and the inclusion map induces the zero map on fundamental groups, which can be seen using the methods from the proof of  \cite[Theorems~5.12 and 5.16]{Roberts_phd}.
In particular this implies that $X^I$ and $LX$ are semilocally simply-connected (a slightly weaker version of this implication follows from a result of Wada \cite{Wada_55}).
This neighbourhood basis has better computational properties for the purposes of this paper than the usual one.
 
If $[f]$ is a 2-track with representative bigon $f\colon I^2 \to X$, let $\ulcorner f\! \urcorner \colon I \to X^I$ be the corresponding path in the mapping space from $s_1[f]$ to $t_1[f]$.

\begin{lemma}\label{lemma:topology_on_2tracks}

  Let $[h] \in \Pi_2(X)_2$ be a 2-track, $U_0,U_1$ basic open neighbourhoods of $s_0[h],t_0[h]\in X$ respectively, and
  \[
		V_0 = N\left(\mathfrak{p}_0,\W^0\right)\quad
		\text{{\rm and }}\quad
		V_1= N\left(\mathfrak{p}_1,\W^1\right)
  \]
  basic open neighbourhoods in $X^I$, where $\W^0 = \{W_i^{(0)}\}_{i=0}^{n_0}$ and $\W^1 = \{W_i^{(1)}\}_{i=0}^{n_1}$. Also assume that  $U_0 \subset W_0^{(0)} \cap W_0^{(1)}$ and  $U_1 \subset W_{n_0}^{(0)} \cap W_{n_1}^{(1)}$ are basic open sets.
  Then the sets
  \begin{multline*}
   \langle [h],U_0,U_1,V_0,V_1 \rangle := \{[f] \in \Pi_2(X)_2\ \mid \ 
   \exists \beta_\epsilon\colon I \to U_\epsilon\text{ {\rm and} } \ulcorner \lambda_\epsilon\! \urcorner \colon I \to V_\epsilon,\;\epsilon=0,1,\\
   \text{{\rm such that }} [f] = [\lambda_1 + (\id_{\beta_1} \bullet (h \bullet \id_{\beta_0})) + \lambda_0]
   \}
  \end{multline*}
  form an open neighbourhood basis for $\Pi_2(X)_2$.

\end{lemma}

\begin{proof}

  We need to show that the axioms for an open neighbourhood basis (See e.g.~conditions a), b) and c') from \cite[Definition~5.6.1]{RBrown_groupoids}) are satisfied. 
  The elements of the putative basic open neighbourhoods $\Y_{[h]} := \langle [h],U_0,U_1,V_0,V_1 \rangle$ are, up to some suppressed bracketing on the whiskering of $[h]$, diagrams of the form
  
\begin{equation*}\label{eq:2-tracks_in_basic_open_nhd}
    \xymatrix{\\
    \\
    x_0 \ar[r]^{\beta_0} \ar@/^4pc/[rrrr]_{\ }="s_2" \ar@/_4pc/[rrrr]^{\ }="t_3"
    &
    s_0[h] \ar@/^1pc/[rr]^{s_1[h]}="t_2"_{\ }="s_1" \ar@/_1pc/[rr]_{t_1[h]}="s_3"^{\ }="t_1"
    &&
    t_0[h] \ar[r]^{\beta_1} & x_1\;,\\
    \\
    \\
    \ar@{=>}"s_1";"t_1"^{[h]}
    \ar@{=>}"s_2";"t_2"_{[\lambda_0]}
    \ar@{=>}"s_3";"t_3"_{[\lambda_1]}
    }
\end{equation*}
  in the bigroupoid $\Pi_2(X)$.
  Figure~\ref{fig:basic_opens}\subref{subfig:basic_neighbourhoods_of_Pi_2_2} is a topological viewpoint of the same element of $\Y_{[h]}$, represented again as a cartoon in Figure~\ref{fig:basic_opens}\subref{subfig:schematic_elt_of_basic_nhd_in_2tracks}.

\begin{figure*}
    \centering
    \subfloat[Topological view.]{
        \centering
        \includegraphics[width=0.5\textwidth]{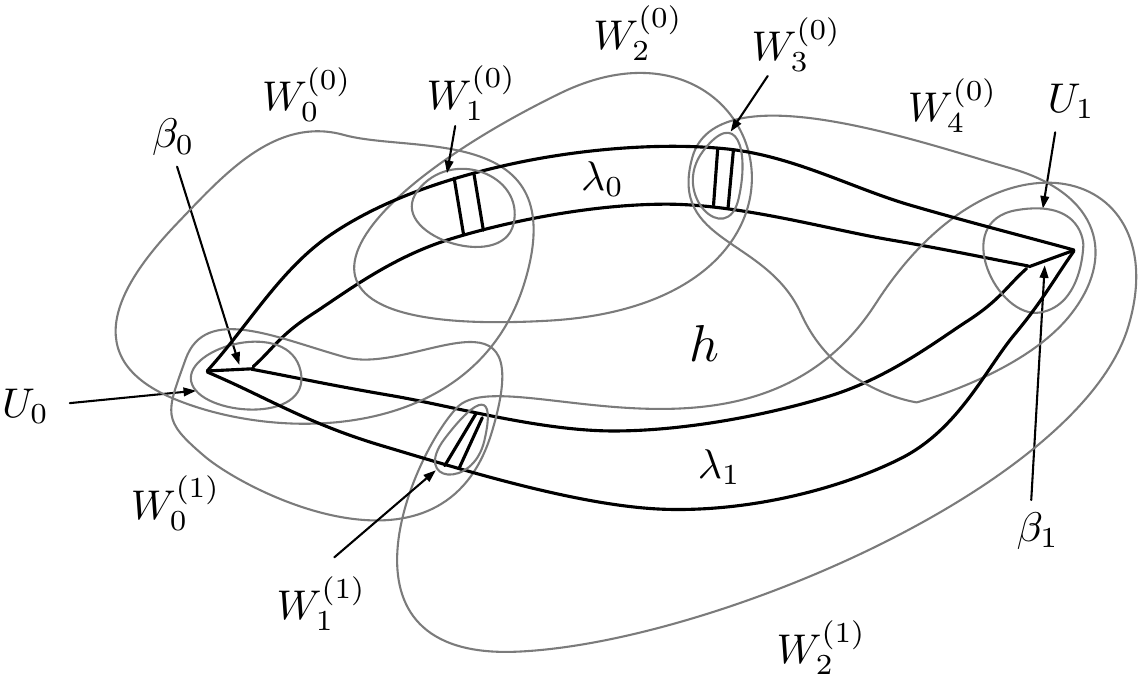}
        \label{subfig:basic_neighbourhoods_of_Pi_2_2}%
    }%
    ~~
    \subfloat[Cartoon view.]{
        \centering
        \includegraphics[width=0.5\textwidth]{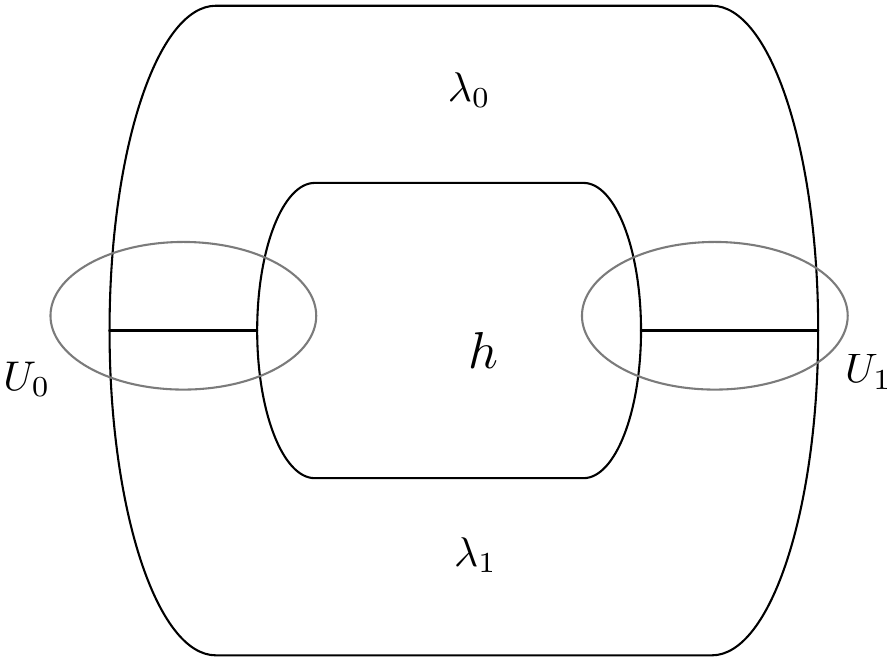}
        \label{subfig:schematic_elt_of_basic_nhd_in_2tracks}%
    }
    \caption{Elements of a basic open neighbourhood in $\Pi_2(X)_2$.}\label{fig:basic_opens}
  \end{figure*}
  
  It is immediate from the definition of $\Y_{[h]}$ that it contains $[h]$. 
  Now assume that $[f]\in \Y_{[h]}$. 
  We need to show that $\Y_{[h]}$ is also a basic open neighbourhood of $[f]$ according to Definition~\ref{def:basis_for_comp_open_topology}. 
  First, notice that $[h]\in\langle [f],U_0,U_1,V_0,V_1 \rangle =: \Y_{[f]}$, since if the 2-track $[f]$ is given by
\[
    \xymatrix{
      \\
      \bullet \ar@/^1pc/[rr]_{\ }="s"^{g_0} \ar@/_1pc/[rr]^{\ }="t"_{g_1} && \bullet
      \ar@{=>}"s";"t"^{[f]}
      \ar@{}[r]|=&
      \bullet \ar[r]^{\beta_0} \ar@/^3.5pc/[rrrr]_{\ }="s_2"^{g_0} \ar@/_3.5pc/[rrrr]^{\ }="t_3"_{g_1}
      &\bullet \ar@/^1pc/[rr]^{s_1[h]}="t_2"_{\ }="s_1" \ar@/_1pc/[rr]_{t_1[h]}="s_3"^{\ }="t_1"
      && \bullet \ar[r]^{\beta_1} & \bullet\\
      \\
      \\
      \ar@{=>}"s_1";"t_1"^{[h]}
      \ar@{=>}"s_2";"t_2"_{[\lambda_0]}
      \ar@{=>}"s_3";"t_3"_{[\lambda_1]}
    }
  \]
  then $[h]$ is given by
  \begin{equation}\label{eq:formula_for_h}
    \raisebox{12ex}{
     \xymatrix@R=4.1ex{
      &&&& \bullet \ar[rr]_{\ }="s_2"^{s_1[h]} && \bullet \ar[dd]_{\beta_1} \ar[rdd]^I_{\ }="s6"
      \\
      \\
      \bullet \ar@/^1pc/[rr]_{\ }="s"^{g_0} \ar@/_1pc/[rr]^{\ }="t"_{g_1} && \bullet
      \ar@{=>}"s";"t"^{[h]}
      \ar@{}[r]|=&
      \bullet \ar[r]="t4_s5"^{\overline{\beta_0}} \ar[ruu]^I_{\ }="s4" \ar[rdd]_I^{\ }="t5" &
      \bullet \ar@/^1pc/[rr]^{g_0}="t_2"_{\ }="s_1" \ar@/_1pc/[rr]="t6_s7"_{g_1} ="s_3"^{\ }="t_1"
      \ar[uu]_{\beta_0} \ar[dd]^{\beta_0}
      && \bullet \ar[r]_{\overline{\beta_1}} & \bullet
      \\
      \\
      &&&& \bullet \ar[rr]^{\ }="t_3"_{t_1[h]} &&\bullet \ar[uu]^{\beta_1} \ar[uur]_I^{\ }="t7" \\
      \ar@{=>}"s_1";"t_1"^{[f]}
      \ar@{=>}"s_2";"t_2"_(.4){-[\lambda_0]}
      \ar@{=>}"s_3";"t_3"_(.6){-[\lambda_1]}
      \ar@{=>}"s4";"t4_s5"
      \ar@{=>}"t4_s5";"t5"
      \ar@{=>}"s6";"t6_s7"
      \ar@{=>}"t6_s7";"t7"
    }
    }
  \end{equation}
  where the unmarked 2-arrows are represented by bigons that are contained in $U_0$ or $U_1$, as appropriate.
  We have not shown all the structure morphisms (associators etc.), relying on coherence for bicategories; see for example \cite{Leinster_98}.

  For an arbitrary 2-track $[g]$ in $\Y_{[h]}$, one can substitute the above expression (\ref{eq:formula_for_h}) for $[h]$ in terms of $[f]$ to get that $[g]$ is in the basic open neighbourhood $\Y_{[f]}$. 
  Thus $\Y_{[h]}\subseteq \Y_{[f]}$.
  By symmetry between $[h]$ and $[f]$ we also have $\Y_{[f]}\subseteq \Y_{[h]}$ and the result follows.
  
  The only thing remaining is to show that an intersection
  \begin{equation}\label{top_bigpd_topology1}
    \langle [h],U_0,U_1,V_0,V_1 \rangle \cap \langle [h],U'_0,U'_1,V'_0,V'_1 \rangle
  \end{equation}
  contains a basic open neighbourhood of $[h]$.
  Choose basic open neighbourhoods
  \[
   V''_0 := N\left(\mathfrak{p}_0,\{W_i^{(0)}\}_{i=0}^{n_0}\right)\subset V_0\cap V'_0, \qquad
   V''_1 := N\left(\mathfrak{p}_1,\{W_i^{(1)}\}_{i=0}^{n_1}\right)\subset V_1\cap V'_1
  \]
  in $X^I$ of the paths $s_1[h]$, $t_1[h]$ respectively, and basic open neighbourhoods
  \begin{align*}
   U''_0 \subset U_0\cap U'_0 \cap W_0^{(0)} \cap W_0^{(1)},\\
   U''_1 \subset U_1\cap U'_1 \cap W_{n_0}^{(0)} \cap W_{n_1}^{(1)}\\
  \end{align*}
  in $X$ of the points $s_0[h]$, $t_0[h]$ respectively. 
  The sets $V''_0$, $V''_1$, $U''_0$ and $U''_1$ satisfy the conditions necessary to make the set
  \[
    \langle [h],U''_0,U''_1,V''_0,V''_1 \rangle
  \]
  a basic open neighbourhood of $[h]$ in $\Pi_2(X)_2$. 
  By inspection this is contained in (\ref{top_bigpd_topology1}) as required. We thus have given a topology on $\Pi_2(X)_2$.
\end{proof}

\subsection{Continuity of source and target for the hom-groupoid}
\label{subsec:hom-gpd_source_target_cts}

Now recall that the map $(s_1,t_1)\colon B_2 \to B_1 \times B_1$ for $B$ a bigroupoid factors through $B_1\times_{B_0\times B_0} B_1$. 
In the case of $\Pi_2(X)$, this gives a function
\[
  \Pi_2(X)_2 \to X^I\times_{X\times X} X^I
\]
of the underlying sets.

If $\mathfrak{p}=\{0 < a_1 <\ldots < a_n < 1\}$ and $\mathfrak{q}=\{0 < b_1 < \ldots < b_m < 1\}$  are partitions we introduce the notation $\mathfrak{p}\vee \mathfrak{q}$ for the partition 
\[
  \{0<a_1/2 <\ldots < a_n/2 < (b_1+1)/2 < \ldots < (b_m+1)/2 < 1\}.
\]

\begin{lemma}\label{Pi2_locally_weakly_discrete}

  With the topology from Lemma~\ref{lemma:topology_on_2tracks},
  $(s_1,t_1)\colon \Pi_2(X)_2 \to X^I\times_{X\times X} X^I$ has  open and closed image, and is a covering map of $\im(s_1,t_1)$.  

\end{lemma}

\begin{proof}
	
  Let $\psi\colon LX \to X^I \times_{X\times X} X^I$ be the map $\omega \mapsto \left(\omega\big|_{[0,\frac12]},\overline{\omega\big|_{[\frac12,1]}} \right)$, where we have implicitly identified $[0,\frac12]\simeq I \simeq [\frac12,1]$ by order-preserving homeomorphisms, and as ever $\overline{(\cdot)}$ denotes the reverse path.
  The homeomorphism $\psi$ will be used in what follows to identify loops and pairs of paths with coinciding endpoints.
  If $L_0X$ denotes the (path) component of the null-homotopic loops, then as we are assuming $X$ is semilocally 2-connected, it is locally path connected, and so $\im(s_1,t_1) \simeq L_0X$ is a component of $X^I \times_{X\times X} X^I$. 
  Hence $\im(s_1,t_1)$ is open and closed.

  Recall from Subsection~\ref{subsec:topology_on_2-arrows} that when $X$ is semilocally 2-connected (Definition \ref{def:sl2c}) the space $LX$ (and hence $L_0X$) is semilocally simply-connected, with path-connected basic open neighbourhoods $N^o(\mathfrak{p},\W)$. 
  It is not difficult to see that $(s_1,t_1)$ is an open map, as it sends basic open neighbourhoods in $\Pi_2(X)_2$ to the basic open neighbourhoods of $X^I\times_{X\times X}X^I$ arising from those in Definition~\ref{def:basis_for_comp_open_topology}.
  Let $\omega$ be a point in $L_0X$, corresponding via $\psi$ to the homotopic paths $\gamma_1,\gamma_2\colon I \to X$ from $x$ to $y$. 
  Let $N := N^o(\mathfrak{p},\W)$ be a basic open neighbourhood of $\omega$ in $L_0X$ where
  \begin{align*}
  	\W^1 &= \{W_i\}_{i=1}^n,\\
  	\W^2 &= \{W_i\}_{i=n+2}^k,\ \text{and}\\
  	\W &=  \{W_0\} \sqcup \W^1 \sqcup \{W_{n+1}\} \sqcup \W^2.
  \end{align*}
  (See Definition~\ref{def:basis_for_comp_open_topology} for the conditions the sets $W_i \subset X$ need to satisfy.)
  Without loss of generality we can assume $\mathfrak{p} = \mathfrak{p}_1\vee\mathfrak{p}_2$, such that $N(\mathfrak{p}_1,\W^1)$ and $N(\mathfrak{p}_2,\W^2)$ are basic open neighbourhoods of $\gamma_1$ and $\overline{\gamma_2}$ (the reverse of the path $\gamma_2$) respectively.
  Consider now the pullback
  \[
    \xymatrix{
      N\times_{L_0X} \Pi_2(X)_2 \ar[rr] \ar@<2.5ex>[dd]_{\pi} && \Pi_2(X)_2 \ar[dd]^{(s_1,t_1)}\\
      \\
      \omega\in\; N \ar[rr]^{\subset} && L_0X
    }
  \]
  We want to show there is an isomorphism $N\times_{L_0X} \Pi_2(X)_2 \simeq N\times \Pi_2(X)(\gamma_1,\gamma_2)$, where $\Pi_2(X)(\gamma_1,\gamma_2) := (s_1,t_1)^{-1}(\omega)$.
  For $[h]\in\Pi_2(X)(\gamma_1,\gamma_2)$, define the following basic open neighbourhood:
  \[
    \Y_{[h]} := \langle[h],W_0,W_{n+1},N(\mathfrak{p}_1,\W^1),N(\mathfrak{p}_1,\W^2)\rangle \subset N\times_{L_0X} \Pi_2(X)_2.
  \]
  By definition, the neighbourhoods $N(\mathfrak{p}_1,\W^1)$ and
  $N(\mathfrak{p}_2,\W^2)$ are path-connected, so $\pi_{[h]}\colon \Y_{[h]} \to N$, the restriction of $\pi$ to $\Y_{[h]}$, is surjective. 
  One can also show it is also injective as follows. 

  Let $[k_1],[k_2]\in \Y_{[h]}$ be such that $(s_1,t_1)[k_1] = (s_1,t_1)[k_2]$. 
  We can assume that $k_i$ is in the form $\lambda_1^i + (\id_{\beta_1^i} \bullet (h \bullet \id_{\beta_0^i})) + \lambda_0^i$ as given in Lemma~\ref{lemma:topology_on_2tracks}. 
  Here $\beta_0^i,\beta_1^i$ are paths in $W_0$ and $W_{n+1}$ respectively, with matching endpoints. 
  By the assumption that $X$ is semilocally 2-connected and the definition of the particular basic neighbourhoods from Lemma~\ref{lemma:topology_on_2tracks}, the open sets $W_0,W_{n+1}$ are simply connected, so we can find an endpoint-preserving homotopy $\eta_\epsilon$ from $\beta_\epsilon^1$ to $\beta_\epsilon^2$ for $\epsilon = 0,1$.
  We can then paste these homotopies with $\lambda_0^1$ to get a surface $\Lambda_0$ sharing a boundary with $\lambda_0^2$, corresponding to a pair of paths in $N(\mathfrak{p}_1,\W^1)$ with matching endpoints. 
  Since $N(\mathfrak{p}_1,\W^1)$ is relatively 1-connected, we can find an endpoint-preserving homotopy in $X^I$ between these two paths -- that is, a filler between the surfaces $\Lambda_0$ and $\lambda_0^2$.
  Similarly, we can paste $\eta_0$ and $\eta_1$ with $\lambda_1^1$ to get a surface $\Lambda_1$ sharing a boundary with $\lambda_1^2$; running the argument again, with $N(\mathfrak{p}_2,\W^2)$ gives a filler between the surfaces $\Lambda_1$ and $\lambda_1^2$.
  These two fillers paste together, with the constant homotopy on $h$, to give a boundary-preserving homotopy between $k_1$ and $k_2$, so that $[k_1]=[k_2]$ and $\pi_{[h]}$ is injective.

  Since $(s_1,t_1)$ is an open surjection, $\pi_{[h]}$ is open and hence an isomorphism. 
  Equipping $\Pi_2(\gamma_1,\gamma_2)$ with the discrete topology, we get an induced map 
  \begin{equation}\label{eq:trivialisation_of_2track_cov_space}
   N\times \Pi_2(\gamma_1,\gamma_2) \simeq \coprod_{[h]\in \Pi_2(\gamma_1,\gamma_2)}\Y_{[h]} \to N\times_{L_0X} \Pi_2(X)_2,
   \end{equation} 
  which is an open surjection. 

  The map (\ref{eq:trivialisation_of_2track_cov_space}) is also injective. Since, if $[h], [h'] \in \Pi_2(\gamma_1,\gamma_2)$ and $[g] \in \Y_{[h]}\cap \Y_{[h']}$, then by the proof of Lemma~\ref{lemma:topology_on_2tracks}, $\Y_{[h]}= \Y_{[h']}$. 
  In particular, $[h']\in \Y_{[h]}$, say, so we can run the above argument used for injectivity of $\pi_{[h]}$ again, with $[h]$ and $[h']$, to get that $[h]=[h']$.
  It then follows immediately that $\pi_{[h]}$ and $\pi_{[h']}$ have disjoint images for $[h]\not=[h']$.
  
  Hence (\ref{eq:trivialisation_of_2track_cov_space}) is a bijection and thus a homeomorphism.
  This implies that $\Pi_2(X)_2 \to L_0X$ is a covering space.
\end{proof}

\begin{corollary}
	
	The two composite maps $s_1,t_1\colon \Pi_2(X)_2 \to X^I\times_{X\times X} X^I \xrightarrow{\pr_i} X^I$ are continuous. 

\end{corollary}

Note that $\Pi_2(\gamma_1,\gamma_2) \subset \Pi_2(X)_2$, with the subspace topology, is discrete.
A special case of this is that $\pi_2(X,x)$ is discrete for any choice of basepoint $x$, when given the topology inherited from $\Pi_2(X)$.

\subsection{The hom-groupoid is topological}
\label{subsec:hom_gpd_topological}

\begin{lemma}

  The 2-tracks and paths in a space, with the topologies as above, form a topological groupoid $\cHom{\Pi_2(X)}$ with arrow space $\Pi_2(X)_2$ and object space $X^I$.

\end{lemma}

\begin{proof}
  We have already seen that the source and target maps are continuous. 
  All that is left to show is that the unit map $\id_{(\cdot)}$, composition $+$ and inversion $-(\cdot)$ are continuous. 
  For the unit map, let $\gamma \in X^I$, and $\Y_{\id_\gamma} := \langle \id_\gamma,U_0,U_1,V_0,V_1\rangle$ be a basic open neighbourhood of $\id_\gamma\in \Pi_2(X)_2$. 
  Define $C := \id_{(\cdot)}^{-1}(\Y_{\id_\gamma})$ and consider the image of $C$ under $\id_{(\cdot)}$:
  \begin{align*}
   \id_{(\cdot)}(C) & = \{ \eta \in \Y_{\id_\gamma} \mid \eta = [\lambda_1 + (\id_{\beta_1} \bullet (\id_ \gamma \bullet \id_{\beta_0}) + \lambda_0] = \id_\chi\}\\
   & = \{ \eta \in \Y_{\id_\gamma} \mid \eta = [\lambda_1 +\id_{\beta_1 \bullet
   (\gamma \bullet \beta_0)} + \lambda_0] = [\lambda_1 + \lambda_0] = \id_\chi\}.
  \end{align*}
  Then $s_1(\lambda_1) = t_1 (\lambda_0) = \beta_1 \bullet (\gamma \bullet\beta_0)$, $t_1(\lambda_1) = s_1(\lambda_0) = \chi$ and $\lambda_0 = -\lambda_1 =: \lambda$. 
  As $\lambda_0$ is a path in $V_0$ and $\lambda_1$ a path in $V_1$, we see that $\lambda$ is a path in $V_0 \cap V_1$ which implies $\chi\in V_0\cap V_1$. 
  If we choose a basic neigbourhood $V_2\subset V_0\cap V_1 \subset X^I$ of $\gamma$, then $V_2 \subset \id_{(\cdot)}^{-1}(\Y_{\id_\gamma})$, and so the unit map of $\cHom{\Pi_2(X)}$ is continuous.

  \begin{figure}
  \centering
  \subfloat[Generic 2-track in the image of a basic open neighbourhood under composition.]{%
    \includegraphics[width=0.48\textwidth]{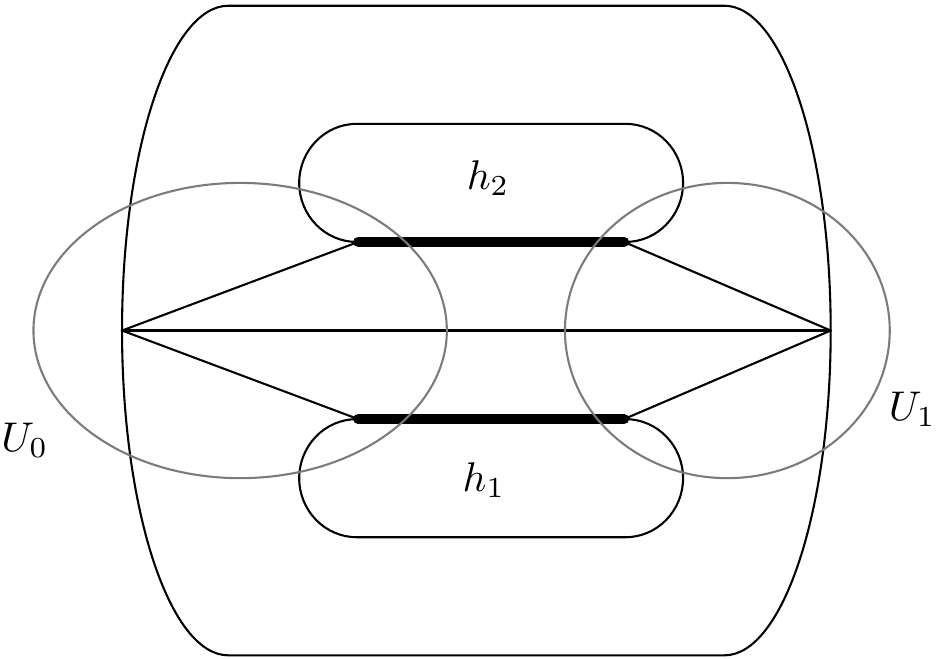}
    \label{subfig:composition_of_2tracks_continuous}
    }
   ~
   \subfloat[Generic 2-track in $\mathcal{I}$.]{%
   \includegraphics[width=0.49\textwidth]{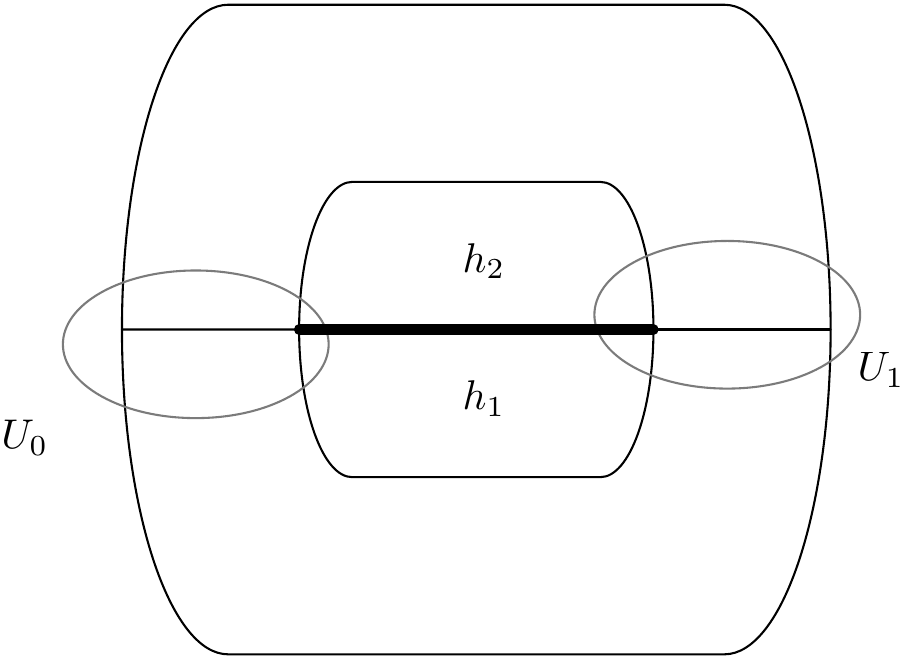}
  \label{subfig:desired_2track_composition}
   }
   \caption{ }\label{fig:vert_comp_cont}
  \end{figure}

  We now need to show the map
  \[
    +\colon \Pi_2(X)_2 \times_{X^I} \Pi_2(X)_2 \to \Pi_2(X)_2
  \]
  is continuous. 
  Let $[h_1]$ and $[h_2]$ be a pair of composable arrows, and consider a basic open neighbourhood $\Y_{[h_2 + h_1]} := \langle [h_2 + h_1],U_0,U_1,V_0,V_2\rangle$.
  Choose a basic open neighbourhood $V_1 = N_\gamma(\mathfrak{p},\W)$ of $\gamma = s_1[h_2] = t_1[h_1]$ in $X^I$ such that the open neighbourhoods $U_0$ and $U_1$ are the first and last basic open neighbourhoods in the collection $\W$. 
  Consider the image $\mathcal{I}$ of $\langle [h_2],U_0,U_1,V_1,V_2\rangle\times_{X^I} \langle [h_1],U_0,U_1,V_0,V_1\rangle$ under $+$.
  Figure~\ref{fig:vert_comp_cont}\subref{subfig:composition_of_2tracks_continuous} is a cartoon of what an element in $\mathcal{I}$ looks like.
  The thick lines are identified, and the interiors of the circles are the basic open sets $U_0,U_1 \subset X$. 
  Topologically Figure~\ref{fig:vert_comp_cont}\subref{subfig:composition_of_2tracks_continuous} is a disk with a cylinder $I\times S^1$ glued to it along some $I\times \{\theta\}$. 
  For this 2-track to be an element of our original neighbourhood $\Y_{[h_2 + h_1]}$ we need to show that the surface that goes `under' the cylinder is homotopic (with fixed boundary) to the one that goes `over' the cylinder, i.e. that there is a filler for the cylinder. 
  Then a generic 2-track $[f_2+f_1]\in \mathcal{I}$ is equal to one of the form
  \[
    [\lambda_1 + (\id_{\beta_1} \bullet ((h_2+h_1) \bullet \id_{\beta_0})) + 
    \lambda_0]\in \Y_{[h_2 + h_1]},
  \]
  which is pictured in Figure~\ref{fig:vert_comp_cont}\subref{subfig:desired_2track_composition}.
  The trapezoidal regions in Figure~\ref{fig:vert_comp_cont}\subref{subfig:composition_of_2tracks_continuous} correspond to paths in $V_1$, which under the identification of the marked edges paste to form a loop in $V_1\subset X^I$. 
  As $X^I$ is semilocally 1-connected, there is a filler for this loop in $X^I$. 
  This implies that there is the homotopy we require, and so composition in $\cHom{\Pi_2(X)}$ is continuous.

  It is clear from the definition of the basic open neighbourhoods of $\Pi_2(X)_2$ that the image of the neighbourhood $\langle [h],U_0,U_1,V_0,V_1 \rangle$ under inversion is $\langle [-h],U_0,U_1,V_1,V_0 \rangle$,
  and so inversion is manifestly continuous.
\end{proof}

\subsection{The fundamental bigroupoid is topological}
\label{subsec:top_bigpd}

The maps $\ev_0,\ev_1\colon X^I \to X$ give us a functor $\cHom{\Pi_2(X)} \to \disc(X\times X)$ of topological groupoids. 
We now have all the ingredients for a topological bigroupoid, but first a lemma about pasting open neighbourhoods of paths with matching endpoints.

Let $\gamma_1,\gamma_2\in X^I$ be paths such that $\gamma_1(1) = \gamma_2(0)$ and let $N_1 := N_{\gamma_1}(\mathfrak{p}_1,\W^1)$, $N_2 := N_{\gamma_2}(\mathfrak{p}_2,\W^2)$ be basic open neighbourhoods. 
For an open set $U \subset W^1_n\cap W^2_m$ (these being the last open sets in their respective collections), define subsets of $X^I$,
\[
  M_1 := \{\eta \in N_1\mid \eta(1) \in U\},\quad M_2 := \{\eta \in N_2\mid \eta(0) \in U\}.
\]
We define the pullback $M_1\times_X M_2$ as a subset of $X^I \times_X X^I$, where this latter pullback is by the maps $\ev_0, \ev_1$. 
The proof of the following lemma should be clear.

\begin{lemma}\label{join_of_nhds_in_XI}

  The image of the set $M_1 \times_X M_2$ under concatenation of paths is the basic open neighbourhood
  \[
	N_{\gamma_2\cdot\gamma_1}(\mathfrak{p}_1\vee\mathfrak{p}_2,\W^1\sqcup \{U\} \sqcup \W^2) \subset X^I.
  \]

\end{lemma}

\noindent We shall denote the image of $M_1\times_X M_2$ as in the lemma by $M_1 \#_U M_2$.

\begin{theorem}\label{thm:we_have_a_top_bigpd}
  
  $\Pi_2(X)$ is a topological bigroupoid.

\end{theorem}

\begin{proof}
  We need to show that the identity assigning functor
  \[
    X \to \cHom{\Pi_2(X)},
  \]
  the concatenation and reverse functors,
  \begin{align*}
   \bullet&\colon\cHom{\Pi_2(X)}\times_{X} \cHom{\Pi_2(X)} \to \cHom{\Pi_2(X)},\\
   \overline{(\cdot)}&\colon \cHom{\Pi_2(X)} \to \cHom{\Pi_2(X)},
  \end{align*}
  and the structure maps in (\ref{structure_maps}) are continuous. 
  In showing these functors are continuous, the only part that needs careful attention is the continuity of the arrow component of the concatenation functor; the rest follows from standard results about path spaces.

  Let $\langle [h_2\bullet h_1],U_0,U_1,V_0,V_1 \rangle$ be a basic open
  neighbourhood in $\Pi_2(X)_2$, where we have the basic open neighbourhoods
  \[
    V_0 = N\left(\mathfrak{p}_0,\W^0\right),\quad V_1 = N\left(\mathfrak{p}_1,\W^1\right)
  \]
  of $s_1[h_2\bullet h_1]$ and $t_1[h_2\bullet h_1]$ in $X^I$ where
  \[
    \W^0 = \{W_i^0\}_{i=0}^n, \qquad \W^1 = \{W_j^1\}_{j=0}^m, \quad n,m \geq 3.
  \]
  We can assume that $\mathfrak{p}_0 = \mathfrak{q}^0_1 \vee \mathfrak{q}^0_2$ and $\mathfrak{p}_0 = \mathfrak{q}^1_1 \vee \mathfrak{q}^1_2$ for partitions $\mathfrak{q}^\epsilon_i$ given as follows:
  \begin{align*}
   \mathfrak{q}^0_1&\colon\{a_1,\ldots,a_k\},\\
   \mathfrak{q}^0_2&\colon\{a_{k+2},\ldots,a_n\},\\
   \mathfrak{q}^1_1&\colon\{b_1,\ldots,b_l\},\\
   \mathfrak{q}^1_2&\colon\{b_{l+2},\ldots,b_m\}.
  \end{align*}
  We now define the neighbourhoods
  \begin{align*}
   V^0_1 & := N\left(\mathfrak{q}^0_1,\{W^0_i\}_{i=0}^k\right),\qquad
   V^0_2 := N\left(\mathfrak{q}^0_2,\{W^0_i\}_{i=k+2}^n\right),\\
   V^1_1 & := N\left(\mathfrak{q}^1_1,\{W^1_j\}_{i=0}^l\right),\qquad
   V^1_2 := N\left(\mathfrak{q}^1_2,\{W^1_j\}_{i=l+2}^m\right).
  \end{align*}
  of $s_1[h_1], s_1[h_2]$ (first row), $t_1[h_1]$ and $t_1[h_2]$ (second row), respectively.

  \begin{figure}
  \centering
  \subfloat[Generic 2-track in the image of a basic open neighbourhood under concatenation.]{%
    \includegraphics[width=0.8\textwidth]{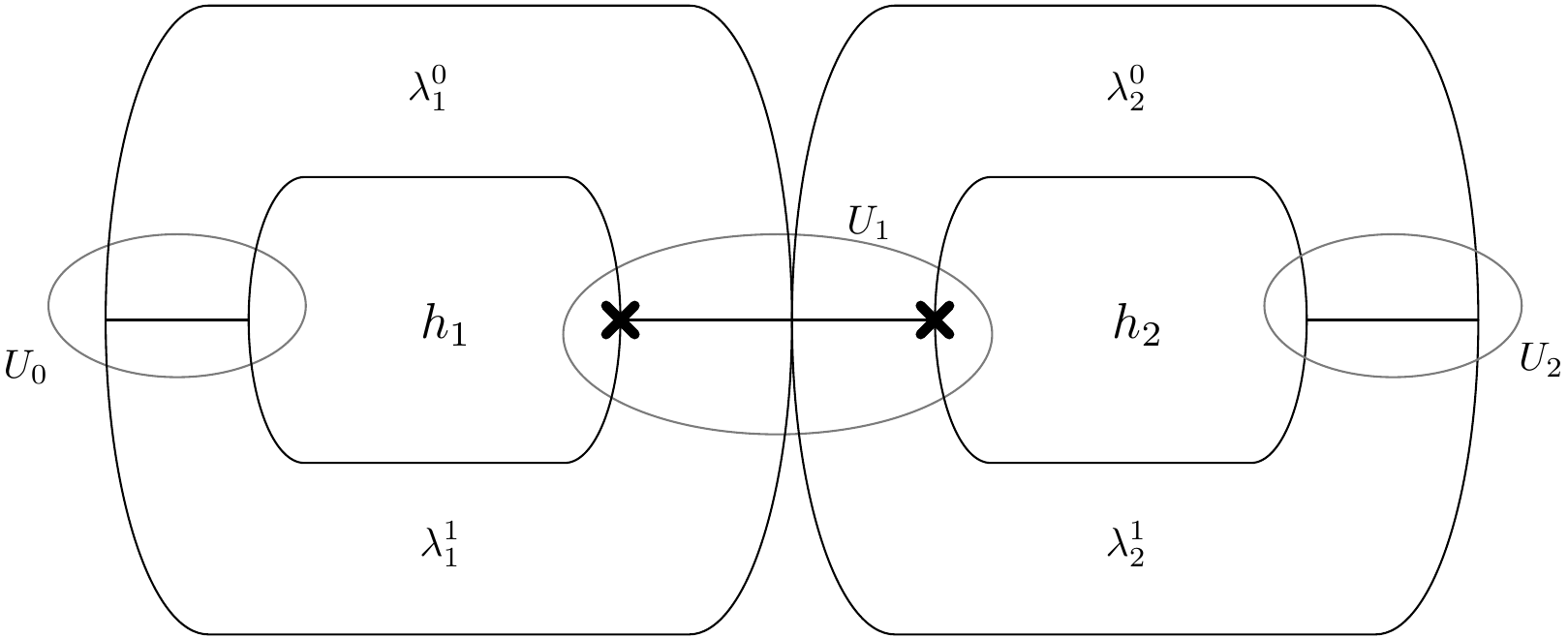}
    \label{subfig:concat_2tracks_cont_a}
  }
  \
  \subfloat[2-track in the image in `standard form'.]{%
    \includegraphics[width=0.65\textwidth]{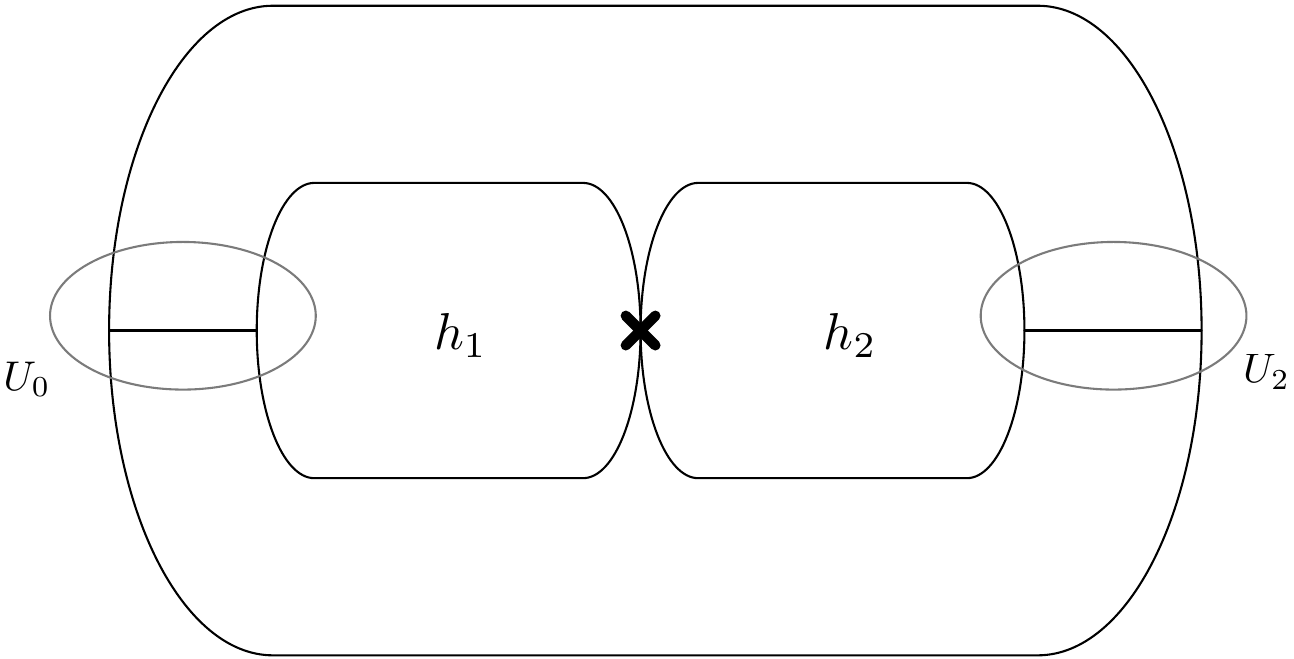}
    \label{subfig:concat_2tracks_cont_b}
  }
  \caption{ }\label{fig:concatenation_of_2tracks_continuous}
  \end{figure}

  Consider the image of the fibred product 
  \[
    \langle [h_1],U_0,U_1,V^0_1,V^1_1
    \rangle\times_X \langle [h_2],U_1,U_2,V^0_2,V^1_2 \rangle
  \]
  under concatenation, any element of which is of the form shown in Figure~\ref{fig:concatenation_of_2tracks_continuous}\subref{subfig:concat_2tracks_cont_a}, where the two points marked with a black cross are identified, so the line between them is a loop in $U_1$. 
  Since the open set $U_1 \subset X$ is 1-connected, there is a filler for this loop, and there is a homotopy between this surface and one of the form showing in Figure~\ref{fig:concatenation_of_2tracks_continuous}\subref{subfig:concat_2tracks_cont_b}.
  Also, by Lemma~\ref{join_of_nhds_in_XI}, the surfaces $\lambda_2^0 \bullet\lambda_1^0$ and $\lambda_2^1 \bullet \lambda_1^1$ are elements of $V^0_1\#_{U_1}V^0_2$ and $V^1_1\#_{U_1} V^1_2$ respectively. 
  Then the image of the open set $\langle [h_1],U_0,U_1,V^0_1,V^1_1 \rangle\times_X \langle [h_2],U_1,U_2,V^0_2,V^1_2 \rangle$ under concatenation is contained in $\langle [h_2\bullet h_1],U_0,U_1,V_0,V_1 \rangle$.

  The assiduous reader will have already noticed that the following relations hold for the (component maps of) the structure morphisms of $\Pi_2(X)$:
  \[
    l = r\circ\overline{(\cdot)},\qquad e = -(i\circ\overline{(\cdot)}).
  \]
  This means that one only needs to check the continuity of $a$ and two of the other four structure maps.

  For the associator $a\colon X^I\times_X X^I\times_X X^I \to \Pi_2(X)_2$, we take a basic open neighbourhood
  \[
    \Y_{a_{\gamma_1\gamma_2\gamma_3}} := \langle a_{\gamma_1\gamma_2\gamma_3},U_0,U_1,V^0,V^1\rangle
  \]
  and by continuity of concatenation of paths choose a basic open neighbourhood $N$ of $(\gamma_1,\gamma_2,\gamma_3)$ in $X^I\times_X X^I\times_X X^I $ whose image under the composite
  \[
    X^I\times_X X^I\times_X X^I \xrightarrow{a}\Pi_2(X)_2 \xrightarrow{(s_1,t_1)} X^I\times_{X\times X}X^I
  \]
  is contained in $V^0\times_{X\times X}V^1$. Also let $U \subset X^I \times_X X^I \times_X X^I$ be a basic open neighbourhood whose image under
  \[
    X^I \times_X X^I \times_X X^I \xrightarrow{a}\Pi_2(X)_2 \xrightarrow{(s_1,t_1)}
    X^I\times_{X\times X}X^I \to X \times X
  \]
  is contained within $U_0\times U_1$. Then if $N' \subset N \cap U$ is a basic open neighbourhood of $(\gamma_1,\gamma_2,\gamma_3)$, its image under $a$ is contained in $\Y_{a_{\gamma_1\gamma_2\gamma_3}}$, so $a$ is continuous.

  The continuity of the other structure maps is proved similarly, and left as an exercise for the reader.
\end{proof}

It is expected that for a reasonable definition\footnote{One possible approach---too much of a diversion to consider here---is to consider the projective local model structure on simplicial sheaves on $\Top$, and the restriction of this to the subcategory of (sheaves represented by nerves of) topological bigroupoids.}  of a weak equivalence of bicategories internal to $\Top$, the canonical 2-functor $\Pi_2^\delta(X) \to \Pi_2(X)$, where recall that $\Pi_2^\delta(X)$ is equipped with the discrete topology, is such a weak equivalence. 
In any case, we can define strict 2-functors between topological bigroupoids, and these are the only such morphisms we shall need here.

\begin{theorem}

  There is a functor
  \[
    \Pi_2\colon \Top_{sl2c} \to \Bigpd(\Top)
  \]
  given on objects by the construction described above, which lifts the fundamental bigroupoid functor $\Pi_2$ of Stevenson and Hardie--Kamps--Kieboom.

\end{theorem}

\begin{proof}
  We only need to check that the strict 2-functor $f_*\colon \Pi_2(X) \to \Pi_2(Y)$ induced by a map $f\colon X\to Y$ in continuous. 
  Recall from \cite{HKK_01} that this strict 2-functor is given by $f$ on objects and post-composition with $f$ on 1- and 2-arrows. 
  We then just need to check that this is continuous on 2-arrows, as it is obvious that it is continuous on objects and 1-arrows.

  Let $\Y_{[f\circ h]} := \langle [f\circ h],U_0^Y,U_1^Y,V_0,V_1\rangle$ be a basic open neighbourhood in $\Pi_2(Y)_2$, and choose basic open neighbourhoods $W_\epsilon \in f^{-1}(V_\epsilon)$ in $X^I$ for $\epsilon = 0,1$. 
  If $\W_0 = \{\coprod_{i=0}^n W_i^{0}\}$ and $\W_1 = \{\coprod_{i=0}^m W_i^{1}\}$, then choose basic open neighbourhoods
  \[
    U_0^X \subset f^{-1}(U_0^Y) \cap W_0^0\cap W_0^1,\quad U_1^X \subset f^{-1}(U_1^Y) \cap W_n^0\cap W_m^1
  \]
  in $X$. 
  It is then clear that $f_*(\langle [h],U_0^X,U_1^X,W_0,W_1\rangle) \subset \Y_{[f\circ h]}$, and so $f_*$ is a continuous 2-functor.
\end{proof}

B\'enabou described in \cite{Benabou} a functor $\Bigpd \to \Gpd$ sending a bigroupoid to the groupoid with the same objects, and isomorphism classes of 1-arrows for arrows. 
Since $\Top$ is cocomplete we can perform the same construction for topological bigroupoids, to get a functor $\Bigpd(\Top) \to \Gpd(\Top)$.

\begin{corollary}

	The composite $\Top_{sl2c} \xrightarrow{\Pi_2} \Bigpd(\Top) \to \Gpd(\Top)$ coincides with the topological fundamental groupoid functor of Brown--Danesh-Naruie.

\end{corollary}

\section{Local triviality}
\label{sec:local_triv}

Recall that for a topological groupoid $\Gamma$ the \emph{source fibre} at an object $p\in \Gamma_0$ is the space $s^{-1}(p)\subset \Gamma_1$.
It follows that the topological group $\Aut(p)$ acts freely on $s^{-1}(p)$ and transitively on the fibres of $s^{-1}(p) \to \Gamma_0$.
For a topological bigroupoid $B$, the source fibre at an object $b\in B_0$ is the sub-topological groupoid $S^{-1}(b) \into \cHom{B}$.
The restriction of the functor $T\colon \cHom{B}\to \disc(B_0)$ then makes $S^{-1}(b)$ a topological groupoid over $\disc(B_0)$.

Recall that a topological groupoid $\Gamma_1 \st \Gamma_0$ is \emph{locally trivial} \cite{Ehresmann_59} if for every point $p\in \Gamma_0$ there is an open neighbourhood $U$ of $p$ such that $s^{-1}(p) \to \Gamma_0$ has a local section on $U$.
If $\Gamma$ is transitive and locally trivial, then one gets local sections of $s^{-1}(p) \to \Gamma_0$ around every point of $\Gamma_0$. 
Thus in this case the source fibre is a (locally trivial) principal bundle.

\begin{example}
	For a semilocally simply-connected topological space $X$, the topological groupoid $\Pi_1(X)$ is locally trivial.
	This is equivalent to the fact that one can find local trivialisations of the universal covering space of $X$.
\end{example}

One can then define a notion of local triviality of topological bigroupoids analogous to that of ordinary topological groupoids. 

\begin{definition}\label{def:loc_triv_bigpd}

  Let $B$ be a topological bigroupoid such that $B_0$ is locally path-connected. 
  We say $B$ is \emph{locally trivial} if the following conditions hold:

  \begin{enumerate}[(I)]
    
    \item
      For every point $b\in B_0$ there is an open neighbourhood $U$ of $b$ such that $\Obj(S^{-1}(b)) \to B_0$ has a local section on $U$;

    \item
      The image of $(s_1,t_1)\colon B_2\to B_1\times_{B_0} B_1$ is open and
      closed, and $B_2 \to \im(s_1,t_1)$ admits local sections.

  \end{enumerate}

\end{definition}

As in the case of 1-groupoids, we get local sections around every point of $B_0$ in the case of a transitive (in that all objects are isomorphic) and locally trivial bigroupoid $B$.
This is related to Bakovic's notion of a bigroupoid 2-torsor \cite{Bakovic_phd}.
In \cite[Chapter~5]{Roberts_phd}, locally trivial bigroupoids were shown, in special cases, to give rise to locally trivial \emph{2-bundles}; this is the motivation for the terminology in Definition~\ref{def:loc_triv_bigpd}, together with the analogy of the situation for 1-groupoids.

While Definition~\ref{def:loc_triv_bigpd} seems to be a good analogue of local triviality for bigroupoids, the main example we are dealing with satisfies a stronger condition than \textrm{(II)}.
This is analogous to the case of the topological groupoid $\Pi_1(X)$, which has the property that $(s,t)\colon \Pi_1(X) \to X\times X$ has discrete fibres.

\begin{definition}\label{def:loc_weakly_discrete}

  A topological bigroupoid $B$ is \emph{locally weakly discrete} if

  \begin{enumerate}
   \item[(II$'$)] The map $(s_1,t_1)\colon B_2\to B_1\times_{B_0\times B_0} B_1$ has discrete (including possibly empty) fibres and is locally trivial.
  \end{enumerate}
  Note that condition (II$'$) implies condition (II) from Definition~\ref{def:loc_triv_bigpd}.

\end{definition}

Recall that a space is locally contractible if it has a neighbourhood basis of contractible open sets. 
We shall call a space \emph{locally relatively contractible} if it has a neighbouhood basis such that the inclusion maps are null-homotopic.

\begin{proposition}

  The bigroupoid $\Pi_2(X)$ is locally weakly discrete, and if $X$ is locally relatively contractible $\Pi_2(X)$ is locally trivial.

\end{proposition}

\begin{proof}
  Lemma~\ref{Pi2_locally_weakly_discrete} shows that $\Pi_2(X)$ is locally weakly discrete, and hence satisfies \textrm{(II)} from Definition~\ref{def:loc_triv_bigpd}. 

  Now assume $X$ is locally relatively contractible. 
  Let $x$ be any point in $X$ and let $U$ be a neighbourhood of $x$ such that $U \into X$ is null-homotopic.
  A homotopy $I\times U \to X$ contracting the inclusion to the base point $x$ then gives a local section $U \to P_xX = \Obj(S^{-1}(x))$, so that $\Pi_2(X)$ satisfies condition \textrm{(I)} and hence is locally trivial.
\end{proof}

\begin{remark}
	As was pointed out by the referee, the singular cubical set of a space can be topologised, and filler operations defined. 
	One may truncate to the level of capturing only 2-dimensional homotopical information, and see what relation this has to the construction of $\Pi_2(X)$.
	Of course, local assumptions on the topology of $X$ are still necessary, otherwise one may not be able to prove the filler operations are continuous.
	The use of the homotopy double groupoid of Brown--Higgins of a stratified space does not immediately appear to be useful for the geometric uses to which the bigroupoid defined above might be put; the truncated and topologised singular cubical set, however, might be useful in defining categorified analogues of covering space constructions as in \cite[Chapter~5]{Roberts_phd}.

  The construction of a homotopy double groupoid that captures the homotopy 2-type of an arbitrary topological space is given in \cite[\S~1.4]{Heredia_phd}.
  This seems like an even more promising strict model that might be lifted to a topological double groupoid and hence approach the geometric constructions for which $\Pi_2$, as given in this paper, was intended.
\end{remark}

\begin{remark}
	Given that the topological fundamental groupoid can be studied in the cases when $\pi_1$ is not discrete, i.e.~when the space at hand is \emph{not} semilocally simply-connected (see for example \cite{Brazas_12}), one wonders whether there is a topological bigroupoid for more general spaces, without the discreteness properties of the current $\Pi_2$.
	Such a structure may not be a topological bigroupoid as defined here, much as one can get a na\"ive topological fundamental group where the multiplication is only \emph{separately} continuous.
	This would give a much finer invariant of spaces and is worth closer consideration.
\end{remark}

\bibliographystyle{amsalpha} 
\bibliography{Pi2inTop_arXiv_final.bbl}

\end{document}